\documentclass[graybox]{svmult}

\usepackage{type1cm}        % activate if the above 3 fonts are
\usepackage{makeidx}         % allows index generation
\usepackage{graphicx}        % standard LaTeX graphics tool
\usepackage{multicol}        % used for the two-column index
\usepackage[bottom]{footmisc}% places footnotes at page bottom

\makeindex             % used for the subject index
\usepackage{amsmath,amsfonts}
\usepackage{array}
\usepackage[caption=false,font=normalsize,labelfont=sf,textfont=sf]{subfig}
\usepackage{textcomp}
\usepackage{stfloats}
\usepackage{url}
\usepackage{verbatim}
\usepackage{graphicx}
\usepackage{cite}
\usepackage{amsbsy}
\usepackage{amsmath}
\usepackage{amssymb}
\usepackage{bm}
\usepackage{soul}
\usepackage{lipsum}% just to generate some text
\def\M3AS{Math.\ Models\ Methods\ Appl.\ Sci.}

\newtheorem{thm}{Theorem}

\newcommand{\prf}{\noindent{{\bf Proof} :\ }}
\newcommand{\QED}{\vrule height 1.4ex width 1.0ex depth -.1ex\ \medskip}
\newcommand{\ls}{\textcolor{red}}
\begin{document}

\title*{Robust path following for autonomous vehicles with spatial PH quintic splines}
\titlerunning{Path Following for autonomous vehicles with PH splines}
\author{Vincenzo Calabr\`o~\orcidID{0000-0003-1246-3128}\\ Carlotta Giannelli\orcidID{0000-0002-5137-1405}\\ Lorenzo Sacco\orcidID{0000-0003-2458-5421}
\\Alessandra Sestini\orcidID{0000-0003-4871-1453}}
\institute{Vincenzo Calabr\`o~ \at MDM Team S.r.l, Florence, Italy \email{vincenzo.calabro@mdmteam.eu}
\and C. Giannelli, L. Sacco, and A. Sestini \at Dipartimento di Matematica e Informatica ``Ulisse Dini'', Universit\`a degli Studi di Firenze, Florence, Italy \email{carlotta.giannelli@unifi.it, lorenzo.sacco@unifi.it, alessandra.sestini@unifi.it}}

\maketitle

\abstract{The distinctive feature of a polynomial parametric speed let polynomial Pythagorean-hodograph (PH) curves be attractive for the design of accurate and efficient application algorithms.
We propose a robust path following scheme for the construction of smooth spatial motions by exploiting PH spline curves. In order to cover a general configuration setting, we present a guidance law which is suitable both for fully-actuated and (more common) under-actuated vehicles, which cannot control all the degrees of freedom. The robustness of the guidance law is enhanced by also taking into account the influence of wind or currents into the equations of motion. A selection of numerical experiments validates the effectiveness of the  control strategy when $C^1$ spatial PH quintic interpolants are suitably considered for both kinematic and dynamic simulations.}

% Comment out for final accepted paper submission
%\linenumbers

%%%%%%%%%%%%%%%%%%%%%%%%%%%%%%%%%%%%%%%%%%%%%%%%%%%%%%%%%%%%%%%%%%%%%
% Introduction
\section{Introduction}
\label{sec:01}
In different application frameworks, from the design of suitable paths for unmanned aerial vehicles (UAVs) or autonomous underwater vehicles (AUVs) to algorithms for computer numerical control machines, a path planning strategy should provide an optimal trade-off between the accuracy of the prescribed trajectory, the efficiency of the control algorithms, and the flexibility of the orientation control of a moving rigid body. To design efficient techniques, the accuracy and the versatility of the motion are usually penalized by approximation schemes based on simple curve types (e.g., linear or circular segments). This kind of restriction can be overcome by considering Pythagorean-hodograph (PH) curves, see \cite{farouki2008,fgs2019} and references therein, characterized by the distinctive feature of a polynomial parametric speed. Algorithms based on PH curve constructions allow to consider smooth spline paths 
%with different smoothness properties, while simultaneously guaranteeing an accurate and efficient arc length computation 
suitable for real-time applications and robust arrival time estimations. Efficient techniques for time computation are crucial in several applications, as the coordinated arrival of vehicle swarms \cite{fgms2018,gkap2009} and the logistical management of different kind of vehicles \cite{znl2015,psb2021}.

In this work we focus on marine robotics, see e.g. \cite{fossen2011}. The model typically used to simulate submerged AUVs kinematics has six degrees of freedom and is generally associated with two ordinary  differential equations (ODEs). First, it is necessary to define two reference frames: the \textit{navigation} reference frame, an inertial frame with a fixed origin, and the \textit{body-fixed} reference frame, whose origin is fixed on a point of the vehicle, while the axes are integral with its movements. By using the quaternion algebra\footnote{In this paper we follow the quaternion notation considered in \cite{farouki2008,fgs2019}.},
%(see Appendix~\ref{sec:appendix} for a brief overview of the basic notations)
the first equation of the considered kinematic model is an ordinary differential equation (ODE) in $\mathbb{R}^3$ which expresses the change $ \dot{\bm{\eta}}$ of body position in the navigation reference system \textcolor{black}{through} a rotation defined by the unit quaternion $\mathcal{Q} = \mathcal{Q}(t)$ as
\begin{equation}
\label{fossen_kinematics_pos} 
\dot{\bm{\eta}} =  \mathcal{Q}\, \bm{v}^b_r \,\mathcal{Q}^* + \dot{\bm{\eta}}_c\,, 
\end{equation}
where $\bm{v}^b_r = \left( u_r^b,v_r^b,w_r^b\right)^\top$ is the velocity of the vehicle with respect to the body-fixed reference frame and $\dot{\bm{\eta}}_c \in \mathbb{R}^3$ is a constant drift velocity of the fluid \cite{fossen2012}. 
%Note that any spatial rotation that maps vectors $\bm{v}^b_r \in \mathbb{R}^3$ expressed in the body-fixed reference system to the navigation reference system can be compactly written as
%$\mathcal{Q}\,\bm{v}^b_r\, \mathcal{Q}^*$ in terms of a unit quaternion $\mathcal{Q} \in \mathbb{H}$ which represents the orientation of the vehicle.
The second equation of the model is an ODE in the quaternion algebra specifying the change of body orientation, 

\begin{equation} 
\label{fossen_kinematics_ori} 
\dot{\mathcal{Q}} =   \frac{1}{2}\, \bm{\omega}\, \mathcal{Q} =   \frac{1}{2}\,   \mathcal{Q}\, \bm{\omega}^b \,, 
\end{equation}
where $\bm{\omega}$ is the angular velocity vector of the body in the navigation reference system and $\bm{\omega}^b =  \mathcal{Q}^*\bm{\omega}\mathcal{Q}$ gives its representation in the body-fixed frame. 

Guidance, navigation, and control are three fundamental parts in the design of control systems for autonomous or unmanned vehicles. In particular, the guidance module is responsible of providing kinematics reference to be followed. 
%The navigation module instead is responsible of estimating the kinematic state of the vehicle (geodetic location, orientation and high order differential states). Finally, the control module combines the results of the previous ones and performs the necessary calculations to allocate forces using available effectors (thrusters, wings, rudder, etc.). Since we do not address navigation and control aspects 
In this paper, we focus on the guidance module \cite{bf2009}, considering a  \emph{path following} approach, see, e.g.,  \cite{llj2008}. 
%In trajectory tracking, a software module is responsible of generating a time-based desired trajectory profile including position $\bm{\eta}_d$, orientation $\mathcal{Q}_d$, angular velocity $\bm{\omega}^b_d$, and high order kinematics states. The trajectory tracking problem requires that the vehicle must \emph{track} both time and kinematics states: the vehicle should be at a certain time in a certain configuration (position/orientation). Modern control solutions, such as dynamic positioning systems for surface vessels and underwater vehicles, are using this approach including additional kinematics constraints for velocity, acceleration, and jerk (acceleration time derivative) \cite{cal2019}. 
%The trajectory tracking strategy is outlined in the block diagram on the left of Fig.~\ref{figure:ttvspf}. 
In this case a software module is responsible of generating a suitable velocity profile $\dot{\bm{\eta}}_d$ to be followed so that the vehicle moves along a desired geometric path without any particular time constraint. 
As our goal is to use the autonomous vehicles in presence of unknown external disturbance, we aim at choosing path following for its intrinsic capability of dealing with these scenarios. We present a guidance law which is suitable both for fully-actuated and (more common) under-actuated vehicles, which cannot control all the degrees of freedom. 
%The robustness of the guidance law is enhanced by also taking into account the influence of wind or currents into the equations of motion. Periodic and systematic updates of estimated time of arrival, can be provided by the system to the operator, to reduce the time-uncertainty of the operation. In order to make these time estimates computationally efficient, PH curves are suitably exploited, granting the capability of computing curve length without involving numerical integration techniques.
%

\textcolor{black}{The structure of the paper is as follows. Section~\ref{sec:02} briefly introduces Pythagorean-hodograph spline curves, while Section~\ref{sec:03} and Section~\ref{sec:04} present the basic and the extended versions of the guidance law, respectively. Kinematic and dynamic simulations are then introduced and commented in Section~\ref{sec:05}. Finally, Section~\ref{sec:07} concludes the paper.}

%%%%%%%%%%%%%%%%%%%%%%%%%%%%%%%%%%%%%%%%%%%%%%%%%%%%%%%%%%%%%%%%%%%%%
% PH
\section{Spatial PH interpolating spline paths}\label{sec:02}

We want to construct a $C^1$ smooth spatial spline path $\bm{\eta}_p(u)$ that interpolates a given ordered set of points $\bm{p}_k \in \mathbb{E}^3, k=0,\ldots,N$,   and corresponding   tangent vectors $\bm{d}_0,\ldots,\bm{d}_N$  (first order Hermite interpolation) so that
\[
\bm{\eta}_p(u_k) = \bm{p}_k, \qquad
\bm{\eta}'_k(u_k) = \bm{d}_k, \qquad
k=0,\ldots,N,
\]
for a certain choice of parameter values $u_0,\ldots,u_N,$ with $u_0 < u_1 < \cdots < u_N.$    Different schemes for associating appropriate interpolation parameters, see, e.g., \cite{farin1997}, and  suitable tangent vectors, see, e.g., \cite{kk2000,acm2001}, to the given points sequence can be considered. Besides these interpolation conditions, we require that each spline segment $ \bm{\eta}_{p,k}$ that defines the spline path,
\[
\bm{\eta}_p(u) := \bm{\eta}_{p,k}(u) = \left(
x_k(u), y_k(u), z_k(u)
\right)^\top
\;\text{if} \; u\in [u_{k-1},u_{k}] \]
for $ k=1,\ldots,N$,
is a polynomial Pythagorean-hodograph (PH) curve, see \cite{farouki2008} and references therein. This means that each hodograph segment $\bm{\eta}'_{p,k}(u) = (x'_k(u),y'_k(u),z'_k(u))^T$ has to satisfy the Pythagorean condition
\[
{x'_k}^{2}(u) + {y'_k}^{2}(u) + {z'_k}^{2}(u) = \sigma_k^2(u),
\]
for some polynomial $\sigma_k$ that represents the curve parametric speed, i.e. $\sigma_k(u)= \Vert \bm{\eta}_{p,k}'(u) \Vert = \mathrm{d}s/\mathrm{d}u$, where $s$ is the cumulative arc length and $ \Vert \cdot \Vert$ denotes the standard Euclidean norm. 

Spatial PH curves admit a compact representation using the algebra of quaternions 
\ls{}%(see  the appendix for a short review) 
and Bernstein polynomials 
$b_i^n(\xi):= \binom{n}{i} \xi^i (1-\xi)^{n-i},\,
\xi\in[0,1],\, 
\text{for } i=0,\ldots,n.$
By focusing on the quintic case, the hodograph segment $\bm{\eta}_{p,k}'$ can be expressed as a quaternion product of the form
\begin{equation}\label{eq:ph}
\bm{\eta}_{p,k}'(u) = \mathcal{A}_k(u)\,\bm{w}\,\mathcal{A}_k^*(u),
\end{equation}
where $\bm{w}$ is any fixed unit vector and $\mathcal{A}_k(u)$ is a quadratic quaternion polynomial which can be expressed in standard B\'ezier form as
\[
\mathcal{A}_k(u) := \sum_{i=0}^2 \mathcal{A}_{k,i} b_i^2\left(\frac{u-u_{k-1}}{h_k}\right)\,, 
\]
with ${\cal A}_{k,i} \in \mathbb{H}$, for $i =0,1,2$ and $h_k := u_k - u_{k-1}$.
By integrating the hodograph \eqref{eq:ph}, we obtain the $k$-th spline segment:
\[
\bm{\eta}_{p,k}(u) = \sum_{i=0}^{5} \bm{q}_{k,i} b_i^5\left(
\frac{u-u_{k-1}}{h_k}\right),
\qquad u\in[u_{k-1},u_k],
\]
defined in terms of the B\'ezier control points 
\begin{align*}
\bm{q}_{k,1} &= \bm{q}_{k,0} + \frac{h_k}{5} 
\mathcal{A}_{k,0}\,\bm{w}\,\mathcal{A}_{k,0}^*,\;
\bm{q}_{k,2}  = \bm{q}_{k,1} + \frac{h_k}{10}
\left(\mathcal{A}_{k,0}\,\bm{w}\,\mathcal{A}_{k,1}^* + 
\mathcal{A}_{k,1}\,\bm{w}\,\mathcal{A}_{k,0}^*\right),\\
\bm{q}_{k,3} & = \bm{q}_{k,2} + \frac{h_k}{30}
\left(\mathcal{A}_{k,0}\,\bm{w}\,\mathcal{A}_{k,2}^*+
4\,\mathcal{A}_{k,1}\,\bm{w}\,\mathcal{A}_{k,1}^*+\right.
\left.\mathcal{A}_{k,2}\,\bm{w}\,\mathcal{A}_{k,0}^*\right),\\
\bm{q}_{k,4} & = \bm{q}_{k,3} + \frac{h_k}{10}
\left(\mathcal{A}_{k,1}\,\bm{w}\,\mathcal{A}_{k,2}^*+
\mathcal{A}_{k,2}\,\bm{w}\,\mathcal{A}_{k,1}^*\right),\;
\bm{q}_{k,5}  = \bm{q}_{k,4} + \frac{h_k}{5}
\mathcal{A}_{k,2}\,\bm{w}\,\mathcal{A}_{k,2}^*\,,
\end{align*}
with $\bm{q}_{k,0}  = \bm{p}_{k-1}$.
In each spline segment, for $u \in [u_{k-1},u_{k}]$, the parametric speed is a quartic polynomial
$
\sigma_k(u) = \mathcal{A}_k(u)\,\mathcal{A}_k(u)^*\,.
$
The distinctive polynomial form of the parametric speed $\sigma$ that characterizes the family of spatial PH curves enables an \emph{exact} and \emph{efficient} computation of the curve arc length, which can be obtained without numerical quadrature by a finite sequence of arithmetic operations on the B\'ezier coefficients. This key feature is strictly connected to the availability of a real-time estimation of the vehicle arrival time. As a matter of fact, it is always possible to compute a first raw estimation of this quantity, starting from the cruise speed and the analytic evaluation of the length of the trajectory. 
%In addition, a polynomial parametric speed is also necessary for a rational unit tangent and, consequently, polynomial spatial curves with rational moving {\it adapted} frames, where one of the frame vectors is the unit tangent to the path
%\footnote{An adapted frame $\bm{e}_1, \bm{e}_2, \bm{e}_3$ along a given path is a triple of mutually orthogonal unit vectors with $\bm{e}_1$ aligned with the unit tangent to the path.}, (as for example the Frenet frame and the rotation-minimizing frame),  are necessarily PH curves. For more details related to theoretical and practical application aspects of PH curves we refer to the recent survey \cite{fgs2019}. 
Many papers recently addressed different application problems with planar or spatial PH curves for the control of autonomous vehicles, see, e.g., \cite{fgms2018} and references therein. In the planar case, PH splines have also been considered in connection with obstacle avoidance techniques \cite{gms2016,dgms2017} and to test a $C^2$ feedrate scheduling algorithm \cite{gms2018}.

We observe that the problem of first order Hermite interpolation by spatial PH quintics admits a two-parameter family of formal solutions \cite{fas2002,sj2005}. 
%In \cite{fgms2008} it was shown that the arc length of the interpolants depends only on one of the two free parameters and alternative criteria to select a suitable interpolant among the complete set of formal solutions were proposed. In \cite{hmk2020} a selection based on extremal interpolants was recently proposed and compared with the so-called CC criterion introduced in \cite{fgms2008}. Since the performances of the two selection schemes are similar but the approach proposed in \cite{hmk2020} generates four interpolants that require additional computations to select the best solution, 
In this paper we consider the CC criterion introduced in \cite{fgms2008} in view of its computational simplicity and of the high shape quality of the interpolants it yields. Note that the approximation order of the interpolation scheme based on this criterion was analyzed in \cite{slm2013}.

% Guidance Law 
\section{Guidance law}\label{sec:03}

As stated in Section \ref{sec:01}, the goal of a path-following scheme is to prescribe the vehicle velocity commands needed to achieve motion control objectives \cite{bf2009}.
This goal is generally achieved by means of a guidance law capable of combining both navigation information, such as vehicle position, with desired path geometrical information, such as tangent direction, parameter value, etc. We will now detail the path parameterization and the base guidance law derived from \cite{bf2005}.

%%%%%%%%%%%%%%%%%%%%%%%%%%%%%%%%%%%%%%%%%%%%%%%%%%%%%%%
\subsection{Path parameterization}
Let us consider a vehicle described in a navigation reference frame. Being the vehicle a rigid body, its position at a certain time instant $t$ can be described by the position $\bm{\eta}(t)= \left( x(t),y(t),z(t) \right)^\top \in \mathbb{E}^3$ of one of its points.
Concerning the geometric path, we adopt the PH spline parametric representation  introduced in the previous section,
$
\bm{\eta}_p(u) = \left(x_p(u),y_p(u),z_p(u)\right)^\top \in \mathbb{E}^3\,, \, u \in [u_0\,,\,u_N]\,.
$
Note that, since only first order derivatives are involved in the calculation of relevant guidance law parts, $\bm{\eta}_p$ is required to be $C^1$ smooth and this motivates our choice of spatial $C^1$ PH spline paths.

\begin{figure}[t]
\centering
\includegraphics[width=0.5\linewidth]{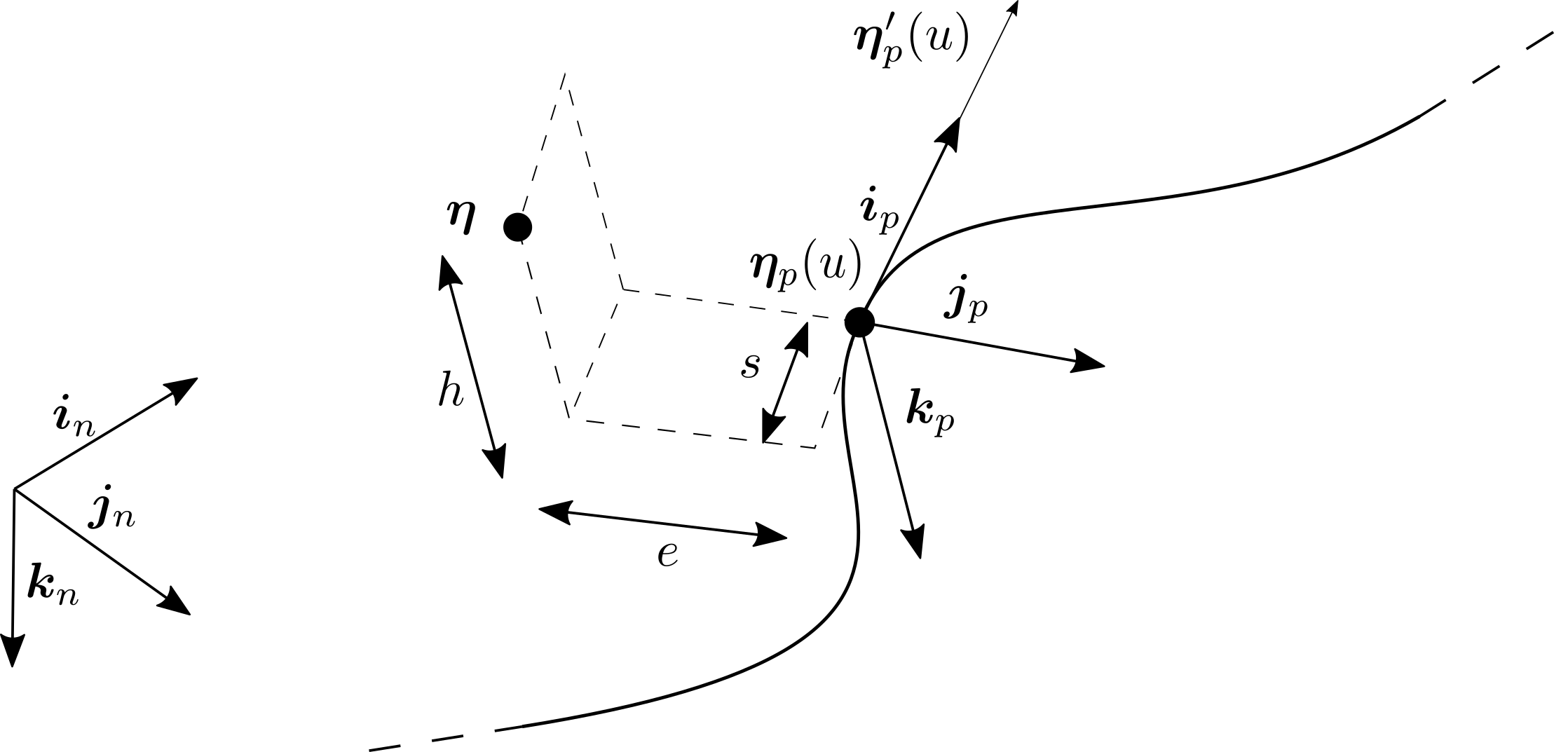}
\caption{Navigation frame (left), path reference frame and track errors $s\,, e\,,h$ (right).}
\label{fig:path_reference_frame}
\end{figure}

We introduce an adapted reference path frame $(\bm{i}_p\,, \bm{j}_p\,, \bm{k}_p)$, with vertex at $\bm{p}_p(u)$ and with the $x$-axis (aligned with $\bm{i}_p$) directed as the tangent to the curve. Since additional requirements on the adapted frame are not considered, the unit quaternion $\mathcal{Q}_p$ which maps vectors represented in the path reference system to the navigation one,
is just required to satisfy the following vector equation (expressed in quaternion algebra):
\begin{equation}
\label{quat_eq}
 \mathcal{Q}_p (1,0,0)^T  \mathcal{Q}_p^* =  \frac{\bm{\eta}_p'}{\Vert \bm{\eta}_p' \Vert}\,,
 \end{equation}
 where standard Euclidean norm is considered. We can then define $\mathcal{Q}_p$ in terms of only two angles (instead of three). In particular, following  \cite{bf2005}, we set
\begin{equation}
\label{unitq}
\mathcal{Q}_p = \mathcal{Q}_z(\chi_p) \mathcal{Q}_y(\nu_p)\,,
\end{equation}
where $\mathcal{Q}_z(\chi_p)$ and $\mathcal{Q}_y(\nu_p)$ are the two unit quaternions 
\begin{equation*}
\mathcal{Q}_z(\chi_p) := \left(\cos \frac{\chi_p}{2}, 0, 0 , \sin \frac{\chi_p}{2} \right)^\top\,, 
\quad 
\mathcal{Q}_y(\nu_p) := \left( \cos \frac{\nu_p}{2}, 0 , \sin \frac{\nu_p}{2}, 0 \right)^\top,
%\mathcal{Q}_z(\chi_p) := \left(\cos \frac{\chi_p}{2}, 0, 0, \sin %\frac{\chi_p}{2} \right)^\top\,, \qquad 
%\mathcal{Q}_y(\nu_p) := \left(\cos \frac{\nu_p}{2}, 0, \sin %\frac{\nu_p}{2}, 0 \right)^\top\,,
\end{equation*}
which define a rotation of an angle $\chi_p$ ({\em azimuth} angle) and $\nu_p$ ({\em elevation} angle) about the $z$ and the $y$ axes, respectively. With the help of quaternion algebra, it can be easily shown that, in order to satisfy (\ref{quat_eq}), the azimuth and elevation angles have to be defined as
\begin{align*}
\chi_p(u) = {\tt atan}_2 \left(x'_p(u)\,,\, y'_p(u) \right),\;
\nu_p(u) = {\tt atan}_2 \left( -z'_p(u)\,,\,\sqrt{x'_p(u)^2 + y'_p(u)^2} \right),
\end{align*}
where ${\tt atan}_2$ denotes the four quadrant inverse tangent. By following the literature in the field, to evaluate the performances of our path following scheme, in our numerical experiments we consider the track error 
\begin{equation} \label{errorvector}
\pmb{\bm{\varepsilon}}^p =
(s,e,h)^\top
:=
\mathcal{Q}_p^*(\bm{\eta} - \bm{\eta}_p(u)){\cal Q}_p \,,
\end{equation}
which is the expression of the error vector $(\bm{\eta} - \bm{\eta}_p(u))$ in the path reference system. The Cartesian components $s\,, e\,, h$ of the track error are usually called {\em along-track} error, {\em cross-track} error, and  {\em vertical-track} error, respectively.
The path reference frame $(\bm{i}_p\,, \bm{j}_p\,, \bm{k}_p)$ and the error components $s\,,e\,,h$ are depicted in Fig. \ref{fig:path_reference_frame}, together with the navigation frame $(\bm{i}_n\,, \bm{j}_n\,, \bm{k}_n)$.

%%%%%%%%%%%%%%%%%%%%%%%%%%%%%%%%%%%%%%%%%%%%%%%%%%%%%%%%%%%%%%%
\subsection{Guidance law}
The objective of the guidance control law is to ensure the convergence of the vehicle to the geometric path. This is done by generating a function $u = u(t)$ and a suitable desired velocity $\dot{\bm{\eta}}_d(t)=\left(\dot{x}_d(t),\dot{y}_d(t),\dot{z}_d(t)\right)^\top \in \mathbb{R}^3\,$ so that
\begin{equation}
\label{def:PathFollowingConvergence}
 \lim_{t \to +\infty} {\Vert\bm{\eta}(t) - \bm{\eta}_p(u(t))\Vert} = 0\,,
\end{equation}
where we assume that the vehicle can have at each time instant the  desired velocity, and consequently $\dot{\bm{\eta}} = \dot{\bm{\eta}}_d.$ We then consider an additional reference frame attached to the vehicle for the development of the guidance law proposed in \cite{bf2005}, which is here revisited.  This frame is called {\it desired reference frame}, since its $x$-axis is in the direction of  the desired velocity. Conversely, the representation of the desired velocity in this frame is $(U_d,0,0)^T$, where $U_d := \Vert \dot{\bm{\eta}}_d \Vert$ is the  desired vehicle speed. A corresponding unit quaternion $\mathcal{Q}_d$ is then necessary to represent the desired velocity in the navigation reference system,
$
\dot{\bm{\eta}}_d := \mathcal{Q}_d
\left(U_d,0,0\right)^\top
\mathcal{Q}_d^*\,.
$
Thus, analogously to $\mathcal{Q}_p$ which was defined through \eqref{unitq}, also $\mathcal{Q}_d$
can be defined as the product of two unit quaternions  $\mathcal{Q}_{d,z}(\chi_d)$ and  $\mathcal{Q}_{d,y}(\nu_d),$ using just two angles $\chi_d$ and $\nu_d$. In \cite{bf2005} it was shown that a suitable choice for these angles is
\begin{align}
\label{eq:DesiredAzimutElevation}
\chi_d = \arctan \left(  \frac{x_{num}}
{x_{den}} \right),
\;
\nu_d = \arcsin(\sin\nu_p \cos\chi_r \cos\nu_r + \cos\nu_p \sin\nu_r),
\end{align}
with
\begin{align*}
x_{num} &= \cos\chi_p \sin\chi_r \cos\nu_r - \sin\chi_p \sin\nu_p \sin\nu_r  + \sin\chi_p \cos\chi_r \cos\nu_p \cos\nu_r,
\\
x_{den} &= -\sin\chi_p \sin\chi_r \cos\nu_r - \cos\chi_p \sin\nu_p \sin\nu_r  + \cos\chi_p \cos\chi_r \cos\nu_p \cos\nu_r,
\end{align*}
\begin{equation} \label{controlangle}
\chi_r(e) = \arctan \left( -\frac{e}{\Delta_e} \right), \qquad \nu_r(h) = \arctan \left( \frac{h}{\Delta_h} \right)\,,
\end{equation}
where $\Delta_h = \mu\sqrt{\Delta_e^2 + e^2}$, while $\Delta_e$ and $\mu$ are two positive guidance parameters. For fixed values of $e$ and $h$ in \eqref{errorvector}, the guidance parameters $\Delta_e$ and $\mu$ influence the direction of convergence along the azimuth and elevation plane respectively. By taking into account (\ref{eq:DesiredAzimutElevation}) and (\ref{controlangle}), let us consider 
$
\dot{u} := (U_d \cos\chi_r \cos\nu_r + \gamma s)/\Vert \bm{\eta}_p' \Vert\,, 
$
where $\gamma$ is another positive control parameter influencing the convergence of the along-track error $s$. For this choice, the Lyapunov theory related to the stability of non linear Cauchy problems ensures that an equilibrium solution $\pmb{\bm{\varepsilon}}^p = \bm{0}$ is obtained for any positive $U_d$ and control parameters \cite{bf2005}. 
In particular, if
\begin{equation}
U_d = \frac{U_0}{\mu\Delta_e}\sqrt{\mu^2(\Delta_e^2 + e^2) + h^2}\,, 
\end{equation} 
where $U_0$ is a given positive constant velocity value, then uniform global exponential stability (UGES) can be achieved, see \cite{tpl2002}. Now, when the vehicle converges to the path, $U_d$ tends to $U_0$. Consequently, the value $U_0$ can be considered a desired speed at steady state. In general, the choice of a suitable value for this steady speed is delegated to the operator.

Note that, in practice, the vehicle will not be able to keep the prescribed velocity, since in the fluid there is a current  that can be modelled as a velocity vector $\dot{\bm{\eta}}_c.$ By denoting with $\dot{\bm{\eta}}_r$ the velocity of the vehicle relative to the fluid, in order to consider a more realistic velocity for the vehicle,  since $\dot{\bm{\eta}} = \dot{\bm{\eta}}_r +\dot{\bm{\eta}}_c$, we should assume $\dot{\bm{\eta}}_r = \dot{\bm{\eta}}_d$. 
Our experiments clearly confirm that this choice does not guarantee the convergence of $\bm{\eta}$ to the path, even under the assumption of a moderate constant current velocity,\footnote{The hypothesis of a constant current is a common assumption in the literature, since it is reasonable during not very long time periods.} see for example the center plots in Fig.~\ref{fig:exm01-01}.
In the next section we will introduce a new {relative desired velocity} $\dot{\bm{\eta}}_{rd}$ to recover the right vehicle behaviour even when a current drifting $\dot{\bm{\eta}}_c$ is considered.
% Extended Gudance Law
% section04.tex

\section{Extended guidance law}\label{sec:04}
In this section we will present an extended version of the basic guidance law to handle two major improvements:
(1) derive a suitable estimation of the ocean current to be effective in realistic scenarios, see Section~\ref{section:Estimation_of_Kinematic_Drifts}, and (2) provide the desired velocity in a suitable form  for general under-actuated underwater vehicles with limited actuation capabilities (i.e., torpedo-like), see Section~\ref{section:map}. These vehicles are typically designed to exploit longitudinal speed at cost of reduced maneuverability. This reflects the limited control possibilities, dealing typically with a single main propeller coupled with a pair of rudders. A simple, yet effective, configuration of this kind reduces the available degrees of freedom of the vehicle to the \emph{surge motion}, used mainly for longitudinal propulsion, the \emph{pitch rotation}, used to steer the vehicle nose up and down to change depth and the \emph{yaw rotation}, used to steer the vehicle nose left and right to change its planar position. For this reason, it is necessary to provide to the controller a desired velocity triplet $( u_{rd}^b, \omega_{yd}^b, \omega_{zd}^b)$, where $u_{rd}^b$ is the surge velocity over the $x$-axis of the body-fixed frame, $\omega_{yd}^b$ is the pitch angular velocity about the $y$-axis of the body-fixed frame and $\omega_{zd}^b$ is the yaw angular velocity about the $z$-axis of the body-fixed frame.

\subsection{Estimation of kinematic drifts}
\label{section:Estimation_of_Kinematic_Drifts}
We will now consider the presence of a current, modelled as a velocity vector of the fluid $\dot{\bm{\eta}}_c$. With this hypothesis, the velocity of the vehicle can not reach the desired value, since it is only possible to assign the relative desired velocity $\dot{\bm{\eta}}_{rd}$ and 
\begin{equation}
\label{eq:relativeProblem}
    \dot{\bm{\eta}} = \dot{\bm{\eta}}_{rd} + \dot{\bm{\eta}}_c,  
\end{equation}
where $\dot{\bm{\eta}}_{rd} = \dot{\bm{\eta}}_{r}$. For this reason, a suitable estimation of the current drift is needed. Let $\dot{\hat{\bm{\eta}}}_c$ be the estimated velocity of the current. By considering the new relative desired velocity as $\dot{\bm{\eta}}_{rd} = \dot{\bm{\eta}}_d - \dot{\hat{\bm{\eta}}}_c$ in (\ref{eq:relativeProblem}), we obtain
\begin{equation}
\label{eq:kinematicsWithEstimate}
\dot{\bm{\eta}} = \dot{\bm{\eta}}_d - \dot{\hat{\bm{\eta}}}_c+ \dot{\bm{\eta}}_c. 
\end{equation}
If the estimate of the current tends to its real value, it is clear that we come back to $\dot{\bm{\eta}} = \dot{\bm{\eta}}_d$. As a consequence, all the convergence results presented in Section \ref{sec:03} still hold true. In the following theorem we derive a suitable estimation for the (constant) current drift that garantees \eqref{def:PathFollowingConvergence}.

\begin{thm}\label{thm:result}
Let the fluid have constant irrotational current $\dot{\bm{\eta}}_c$ and $\dot{\bm{\eta}}$ given by (\ref{eq:kinematicsWithEstimate}). If $\dot{\bm{\eta}}_d$ is computed with the guidance law of Section \ref{sec:03} and
\begin{equation}
\label{currentmodel}
\ddot{\hat{\bm{\eta}}}_c = k_c(\bm{\eta}-\bm{\eta}_p),\text{ with $k_c > 0$}, 
\end{equation}
then 
$
\displaystyle{\lim_{t \to \infty} \Vert \dot{\tilde{\bm{\eta}}}_c \Vert = 0}$ and $\displaystyle{\lim_{t \to \infty} \Vert \pmb{\bm{\varepsilon}}^p \Vert = 0,}
$
where $\dot{\tilde{\bm{\eta}}}_c := \dot{\bm{\eta}}_c - \dot{\hat{\bm{\eta}}}_c$ is the estimation error and $\pmb{\bm{\varepsilon}}^p$ is the track error.
\end{thm}

\prf
The  modified dynamical system representing the new path following strategy is characterized by the state vector $\pmb{\bm{\varepsilon}}^p  $ defined in (\ref{errorvector}) and by $\dot{\tilde{\bm{\eta}}}_c \in \mathbb{R}^3$. We want to demonstrate that the system admits the asymptotically stable equilibrium solution  
$
(\pmb{\bm{\varepsilon}}^p,
\dot{\tilde{\bm{\eta}}}_c
)^\top 
= (\bm{0}, \bm{0} )^\top. 
$
It is trivial to demonstrate that the vanishing of $\pmb{\bm{\varepsilon}}^p$ and $\dot{\tilde{\bm{\eta}}}_c$ leads to $\dot{\pmb{\bm{\varepsilon}}}^p = \bm{0}$ and $\ddot{\tilde{\bm{\eta}}}_c = \bm{0}$, implying that the solution is an equilibrium. For the stability analysis of such equilibrium solution, let us consider the following Lyapunov function:
$
V(\pmb{\bm{\varepsilon}}^p,\dot{\tilde{\bm{\eta}}}_c) := k_c\frac{1}{2}\pmb{\bm{\varepsilon}}^{p\top}\pmb{\bm{\varepsilon}}^p + \frac{1}{2}\dot{\tilde{\bm{\eta}}}_c^\top\dot{\tilde{\bm{\eta}}}_c\,.
$
It holds that $V(\pmb{\bm{\varepsilon}}^p,\dot{\tilde{\bm{\eta}}}_c) \ge 0\,,$ with $V(\pmb{\bm{\varepsilon}}^p,\dot{\tilde{\bm{\eta}}}_c) = 0 \text{ if and only if }\pmb{\bm{\varepsilon}}^p =  \dot{\tilde{\bm{\eta}}}_c = \bm{0}\,.$  
In order to study the behavior of this function, we need now to compute its time derivative
\begin{equation}
\label{eq:lyapunovDerivative}
\dot{V}(\pmb{\bm{\varepsilon}}^p,\dot{\tilde{\bm{\eta}}}_c) = k_c\pmb{\bm{\varepsilon}}^{p\top}\dot{\pmb{\bm{\varepsilon}}}^p + \dot{\tilde{\bm{\eta}}}_c^\top\ddot{\tilde{\bm{\eta}}}_c.
\end{equation}
By considering the definition of $\pmb{\bm{\varepsilon}}^p$ in (\ref{errorvector}), together with (\ref{eq:kinematicsWithEstimate}), we obtain
\begin{align}
\label{eq:errorDerivative}
    \nonumber 
    \dot{\pmb{\bm{\varepsilon}}}^p &= 
    \dot{\mathcal{Q}}_p^*(\bm{\eta} - \bm{\eta}_p(u))\mathcal{Q}_p + \mathcal{Q}_p^*(\bm{\eta} - \bm{\eta}_p(u))\dot{\mathcal{Q}}_p + \mathcal{Q}_p^*(\dot{\bm{\eta}} - \dot{\bm{\eta}}_p(u))\mathcal{Q}_p  
    \\
    & = \dot{\mathcal{Q}}_p^*\, \mathcal{Q}_p\,\pmb{\bm{\varepsilon}}^p + \pmb{\bm{\varepsilon}}^p\mathcal{Q}_p^*\,\dot{\mathcal{Q}}_p + \mathcal{Q}_p^*\,(\dot{\bm{\eta}}_d + \dot{\tilde{\bm{\eta}}}_c - \bm{\eta}_p'\dot{u})\mathcal{Q}_p\,.
\end{align}
 Furthermore, since a constant current implies $\ddot{\bm{\eta}}_c = \bm{0},$ we have that $\ddot{\tilde{\bm{\eta}}}_c = - \ddot{\hat{\bm{\eta}}}_c.$ Replacing the derived expressions of $\dot{\pmb{\bm{\varepsilon}}}^p$ and of $\ddot{\tilde{\bm{\eta}}}_c$ into (\ref{eq:lyapunovDerivative}), we obtain
\begin{align*}
\dot{V}(\pmb{\bm{\varepsilon}}^p,\dot{\tilde{\bm{\eta}}}_c) 
&= k_c \,\pmb{\bm{\varepsilon}}^{p\top}(\dot{\mathcal{Q}}_p^* \mathcal{Q}_p\pmb{\bm{\varepsilon}}^p + \pmb{\bm{\varepsilon}}^p\mathcal{Q}_p^*\dot{\mathcal{Q}}_p + \mathcal{Q}_p^*(\dot{\bm{\eta}}_d + \dot{\tilde{\bm{\eta}}}_c - \bm{\eta}_p'\dot{u})\mathcal{Q}_p) + \dot{\tilde{\bm{\eta}}}_c^\top(- \ddot{\hat{\bm{\eta}}}_c)\\
&= k_c \dot{V}_e(\pmb{\bm{\varepsilon}}^p) + \dot{V}_c(\pmb{\bm{\varepsilon}}^p,\dot{\tilde{\bm{\eta}}}_c), 
\end{align*}
where 
$
\dot{V}_e(\pmb{\bm{\varepsilon}}^p) := \pmb{\bm{\varepsilon}}^{p\top}(\dot{\mathcal{Q}}_p^* \mathcal{Q}_p\pmb{\bm{\varepsilon}}^p + \pmb{\bm{\varepsilon}}^p\mathcal{Q}_p^*\dot{\mathcal{Q}}_p + \mathcal{Q}_p^*(\dot{\bm{\eta}}_d - \bm{\eta}_p'\dot{u})\mathcal{Q}_p)
$
is a term corresponding to the same Lyapunov function used in  \cite{bf2005}, while
$
\dot{V}_c(\pmb{\bm{\varepsilon}}^p,\dot{\tilde{\bm{\eta}}}_c) := k_c\pmb{\bm{\varepsilon}}^{p\top}(\mathcal{Q}_p^* \dot{\tilde{\bm{\eta}}}_c \mathcal{Q}_p) - \dot{\tilde{\bm{\eta}}}_c^\top (\ddot{\hat{\bm{\eta}}}_c)
$
is a new term necessary to study the convergence of the estimation error. For what concern the term $ \dot{V}_e(\pmb{\bm{\varepsilon}}^p)$, expanding all the products and replacing the values of desired velocity $\dot{\bm{\eta}}_c$ and parameter derivative $\dot{u}$, after few calculations we have that
\begin{equation*}
    \dot{V}_e(\pmb{\bm{\varepsilon}}^p) = 
    -\gamma\, s^2 - \frac{U_0}{\Delta_e}e^2 - \frac{U_0}{\mu\,\Delta_e}h^2\,,
\end{equation*}
which is $\leq 0$ (negative semidefinite) since $\dot{\tilde{\bm{\eta}}}_c$ could be non-zero when  $\dot{V}_e(\pmb{\bm{\varepsilon}}^p) = 0$. Let us now verify that the term $\dot{V}_c(\pmb{\bm{\varepsilon}}^p,\dot{\tilde{\bm{\eta}}}_c)$ vanishes. Since $\ddot{\hat{\bm{\eta}}}_c = k_c(\bm{\eta}-\bm{\eta}_p)\,,$ we obtain
$
\dot{V}_c(\pmb{\bm{\varepsilon}}^p,\dot{\tilde{\bm{\eta}}}_c) = k_c (\pmb{\bm{\varepsilon}}^{p\top}(\mathcal{Q}_p^* \dot{\tilde{\bm{\eta}}}_c \mathcal{Q}_p) - \dot{\tilde{\bm{\eta}}}_c^\top (\mathcal{Q}_p \pmb{\bm{\varepsilon}}^p \mathcal{Q}_p^*)).
$
The four vectors involved in this equation represent the same two vectors expressed in the path reference frame as $\pmb{\bm{\varepsilon}}^{p} \, , \, \mathcal{Q}_p^* \dot{\tilde{\bm{\eta}}}_c \mathcal{Q}_p$ and in the navigation reference frame as $ \dot{\tilde{\bm{\eta}}}_c \, , \, \mathcal{Q}_p \pmb{\bm{\varepsilon}}^p \mathcal{Q}_p^*$ respectively. Consequently, $\dot{V}_c(\pmb{\bm{\varepsilon}}^p,\dot{\tilde{\bm{\eta}}}_c) = 0$ and 
$
\dot{V}(\pmb{\bm{\varepsilon}}^p,\dot{\tilde{\bm{\eta}}}_c) =  k_c \dot{V}_e(\pmb{\bm{\varepsilon}}^p) \leq 0\,.
$
We can than conclude that the system is stable. In any case, we need asymptotic stability to have convergence for both current and error estimation. Thus, let us consider the set
$
    \mathcal{E} = \{ (\pmb{\bm{\varepsilon}}^p , \dot{\tilde{\bm{\eta}}}_c) \,\, | \,\, \dot{V}( \pmb{\bm{\varepsilon}}^p , \dot{\tilde{\bm{\eta}}}_c) = 0 \} = 
 \{ (\pmb{\bm{\varepsilon}}^p , \dot{\tilde{\bm{\eta}}}_c) \,\, | \,\, \pmb{\bm{\varepsilon}}^p = 0 \}\,.    
$
We can observe that $\mathcal{E}$ does not contain any trajectory of the system, except the equilibrium solution $\pmb{\bm{\varepsilon}}^p= \dot{\tilde{\bm{\eta}}}_c = \bm{0}$. Indeed, if we had $\pmb{\bm{\varepsilon}}^p = \bm{0}$ but $\dot{\tilde{\bm{\eta}}}_c \neq \bm{0}$, replacing inside (\ref{eq:errorDerivative}) it would result that $\dot{\pmb{\bm{\varepsilon}}}^p \neq \bm{0}$, and so the trajectory would not stay inside $\mathcal{E}$. As a result we can conclude that the system is asymptotically stable by applying the \textit{LaSalle-Krasovskii principle}, see e.g. \cite{ssvg2008}.
\QED
%%%%%%%%%%%%%%%%

\ls{}
%\begin{rmk} The value of the coefficient $k_c$ which appears in the current model introduced in (\ref{currentmodel}) can be arbitrarily chosen, as long as it is positive. Selecting an higher value allows a faster estimation for the error but, as a counter-effect, leads to more evident oscillations and overshoot effects, caused by \eqref{currentmodel} itself, which is commonly denoted as an integral action in control theory. Moreover, once this guidance law is used for a real vehicle, integral actions are able to estimate additional errors, caused by non-modelled dynamical effects. A suitable trade off has then to be considered in the choice of $k_c$.
%\end{rmk}
 
%%%%%%%%%%%%%%%%%%%%
\subsection{Mapping desired guidance into controllable axes}\label{section:map}
The result of Theorem~\ref{thm:result} allows the vehicle to reject errors caused by the presence of a current but the desired relative velocity $\dot{\bm{\eta}}_{rd}$ is a velocity described in the navigation frame. 
In order to map $\dot{\bm{\eta}}_{rd}$ into three body-fixed velocities $u_{rd}^b, \omega_{yd}^b, \omega_{zd}^b$, it is necessary to follow the path with a point different from the center of the body-fixed reference frame. We consider for this purpose the head of the vehicle, whose constant position in the body-fixed reference frame is 
$
\bm{O}^b_h = (
o_x^b,o_y^b,o_z^b
)^\top.
$
We can now compute its position $\bm{\eta}_h$ in the navigation reference frame as
$
\bm{\eta}_h = \bm{\eta} + \mathcal{Q}\, {\bm{O}^b_h}\,\mathcal{Q}^*,
$
where $\bm{\eta}$ is the position of the center of the body-fixed reference frame and $\mathcal{Q}$ is the quaternion representing the attitude of the vehicle. The body-fixed velocity of the head can be computed as
$
\bm{v}_h^b = \mathcal{Q}^*{\dot{\bm{\eta}}_h}\mathcal{Q}= \mathcal{Q}^*(\dot{\bm{\eta}} + \dot{\mathcal{Q}} {\bm{O}^b_h}\mathcal{Q}^* + \mathcal{Q}{\bm{O}^b_h}\dot{\mathcal{Q}}^*)\mathcal{Q}.
$
By substituting the expression for $\mathcal{Q}$ provided by (\ref{fossen_kinematics_ori}) into the previous equation, while  considering that $\mathcal{Q}$ is a unit quaternion and $(\bm{\omega}\mathcal{Q})^* = \mathcal{Q}^* \bm{\omega}^{*} = -\mathcal{Q}^*\bm{\omega} $, we obtain
$
\bm{v}_h^b 
%=\mathcal{Q}^*\,{\dot{\bm{\eta}}}\,\mathcal{Q} + \frac{1}{2}\mathcal{Q}^* \,\bm{\omega}\,\mathcal{Q}\,{\bm{O}^b_h}\,\mathcal{Q}^*\,\mathcal{Q} - \frac{1}{2}\mathcal{Q}^*\,\mathcal{Q} \,{\bm{O}^b_h} \,\mathcal{Q}^*\,\bm{\omega}\mathcal{Q} 
 =\bm{v}^b + \frac{1}{2} \left (\bm{\omega}^b\,{\bm{O}^b_h} - \bm{O}^b_h\,\bm{\omega}^{b} \right)\,,
$
where we recall that $\bm{\omega}^b = \mathcal{Q}^*\bm{\omega}\mathcal{Q}.$ Thus, using the standard cross product, we can also write
$ 
\bm{v}_h^b = \bm{v}^b + \bm{\omega}^b\times\bm{O}^b_h = \bm{v}^b - \bm{O}^b_h\times \bm{\omega}^b\,.
$
 Finally, we can pass to relative velocities:
\begin{equation}
\label{eq:HeadBodyVelocity}
\bm{v}_{rh}^b = \bm{v}_r^b  - {\bm{O}^b_h} \times \bm{\omega}^b,
\end{equation} 
where $\bm{v}_{rh}^b$ is the relative body-fixed velocity of the head and $\bm{v}_r^b$ is the relative body-fixed velocity of the vehicle. In order to deal with an under-actuated vehicle, we can only control the motion along three directions, as mentioned at the beginning of this section. We can then generate the desired velocities in the following form:
\begin{align}
\label{eq:underactuatedVelocities}
 \bm{v}_{rd}^b =
\bm{v}_{r}^b =
\left(
u_{rd}^b,0,0
\right)^\top,
\quad
\bm{\omega}_d^b =
\bm{\omega}^b = 
\left(
0,\omega_{yd}^b,\omega_{zd}^b
\right)^\top.
\end{align}  
By substituting $\bm{v}_{rd}^b $ and $\bm{\omega}_d^b$ into (\ref{eq:HeadBodyVelocity}), after some trivial calculations, we obtain:
\begin{equation}
\label{eq:MapVelocity}
\bm{v}_{rhd}^b =
\bm{v}_{rh}^b 
= \bm{v}_{rd}^b  - {\bm{O}^b_h} \times \bm{\omega}_d^b = \bm{P}({\bm{O}^b_h})
\left(
u_{rd}^b,
\omega_{yd}^b,
\omega_{zd}^b
\right)^\top, 
\end{equation}
with $\bm{P}({\bm{O}^b_h}) = (
1, o_z^b, -o_y^b;
0, 0 , o_x^b;
0, -o_x^b, 0)$.
Since $\bm{O}^b_h$ is the head of the vehicle, $o_x^b \neq 0$ and, consequently, the matrix $\bm{P}({\bm{O}^b_h})$ is invertible. The solution of the linear system in \eqref{eq:MapVelocity} is  then
$
\left(u_{rd}^b,\omega_{yd}^b,\omega_{zd}^b \right)^\top =\bm{P}({\bm{O}^b_h})^{-1}\mathcal{Q}^*(\dot{\bm{\eta}}_d-\dot{\hat{\bm{\eta}}}_c)\mathcal{Q},
$
 where $ \bm{v}_{rhd}^b = \mathcal{Q}^*\dot{\bm{\eta}}_{rd}\mathcal{Q}\,$ and $\dot{\bm{\eta}}_{rd} = (\dot{\bm{\eta}}_d-\dot{\hat{\bm{\eta}}}_c)$. Summarizing and including the current estimation, the equations of the extended guidance law are 
\begin{align*}
&\bm{\eta}_h =  \bm{\eta} + \mathcal{Q}{\bm{O}^b_h}\mathcal{Q}^*, \qquad  \ddot{\hat{\bm{\eta}}}_c = k_c(\bm{\eta}_h-\bm{\eta}_p), \\
&\left(u_{rd}^b,\,\omega_{yd}^b,\,\omega_{zd}^b\right)^\top = \bm{P}({\bm{O}^b_h})^{-1}\mathcal{Q}^*(\dot{\bm{\eta}}_d-\dot{\hat{\bm{\eta}}}_c)\mathcal{Q}\,. 
\end{align*}
%The  block diagram which outlines the key ingredients of both the base and the extended guidance laws is depicted in Fig.~\ref{fig:guidance_fossen}.
%\begin{figure}[t]
%\includegraphics[width=\linewidth]{img/guidance_fossen.png}
%%\includegraphics[scale=.25]{img/guidance_fossen.png}
%\caption{Block diagram of the basic and extended guidance laws.  The initial value for the integration of the parameter $u$ is chosen as the minimum of the parametric interval where the spline path is defined. The initial value for the integration of the estimated current is conventionally chosen as $\dot{\hat{\bm{\eta}}}_c = \bm{0} \, \rm{m}/s$.}
%\label{fig:guidance_fossen}
%\end{figure}
% Numerical Example
% section05.tex
\section{Numerical tests}\label{sec:05}
For the experiments we consider three different scenarios, performing  both kinematic and dynamic simulations to verify the robustness of the scheme. The guidance parameters are set as follows: $\gamma = 1 \,\rm{s}^{-1}$, $U_0 = 0.4 \, \rm{m}/\rm{s}$, $\Delta_e = 5 \, \rm{m}$ and $\mu = 1$ for the base guidance law; $\bm{O}^b_h = (0.8,0,0)^\top \, \rm{m}$ and $\kappa_c=0.015 \,\rm{s}^{-2}$ for the extended guidance law. In all figures related to the paths, the red markers indicate the starting point of the vehicle and the first point of the path. 
When track errors and current estimation are shown, the vertical dashed lines refer to the junction points of the PH spline path.  The black dots indicate the interpolation points. %considered in the construction of the spline path.

\subsection{Kinematic simulations}
\begin{figure*}[!t]
\centering
%\includegraphics[width=.45\linewidth]{img/base_segments_no_current/path_following.pdf}
%\includegraphics[width=.45\linewidth]{img/base_segments_current/path_following.pdf}
%
%\subfigure[guidance law without current]
{\includegraphics[width=.325\linewidth]{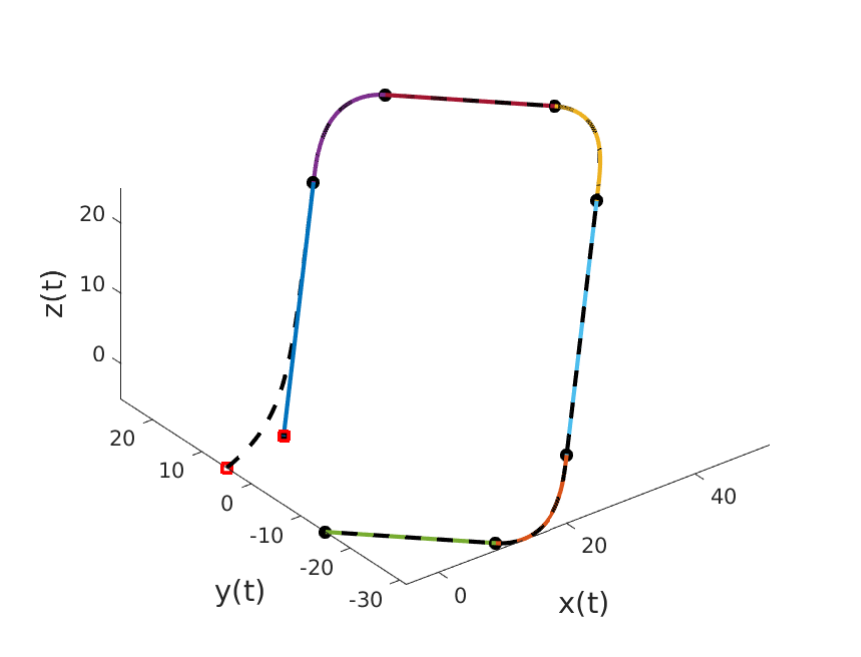}}
%\subfigure[guidance law with current]
{\includegraphics[width=.325\linewidth]{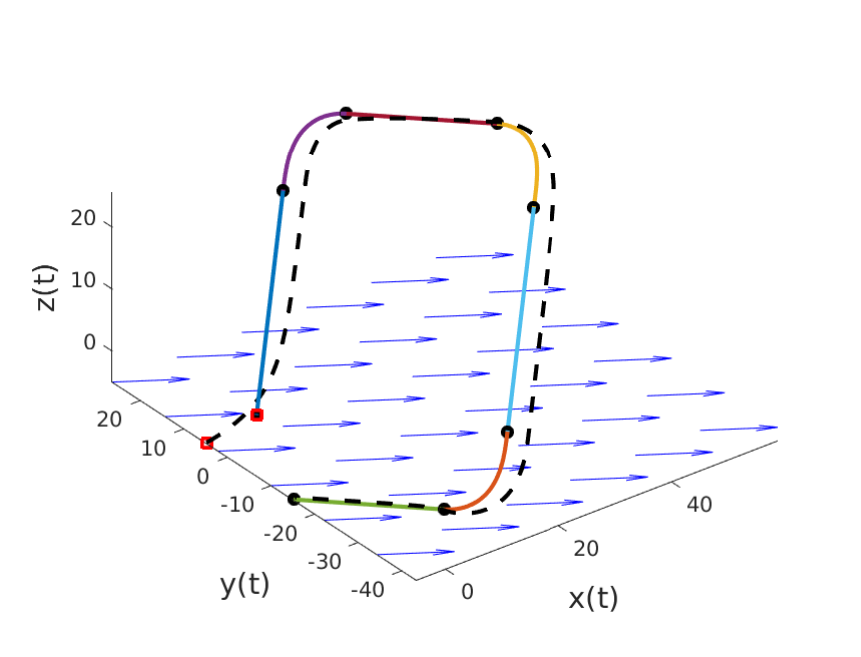}}
%\subfigure[extended guidance law with current]
{\includegraphics[width=.325\linewidth]{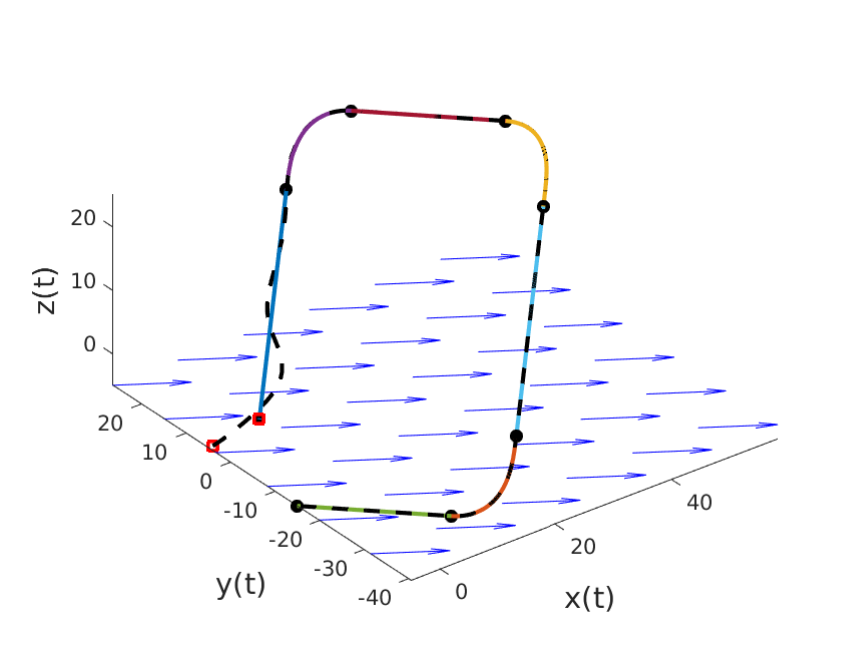}}
\caption{Kinematic simulations: results obtained with the guidance law of Section~\ref{sec:03} and \ref{sec:04} for the first simulation scenario. 
The path to be followed (solid line) is shown together with the path of the vehicle (dashed line) in a simulation without current (left) and with current when the basic guidance law (center) or its extension (right) are considered.}
    \label{fig:exm01-01}
\end{figure*}

\begin{figure}[t]
\centering 
\includegraphics[width=.315\linewidth]{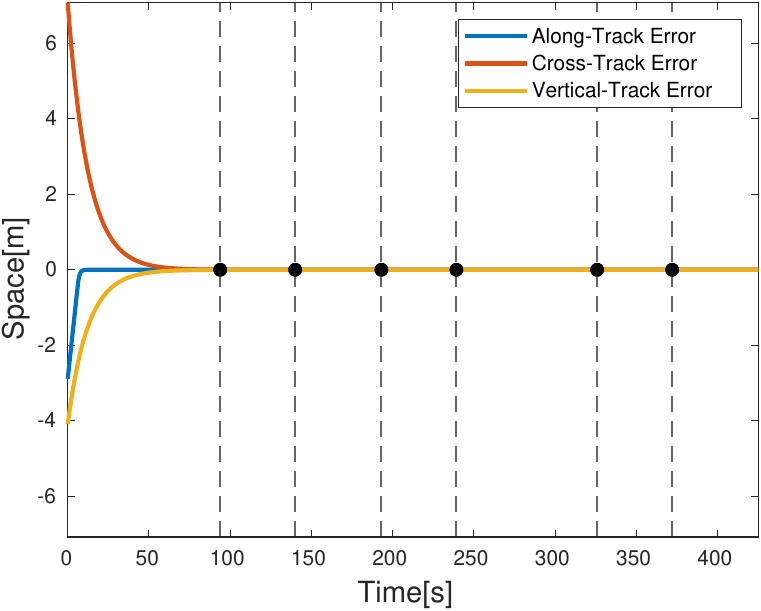}   
\includegraphics[width=.315\linewidth]{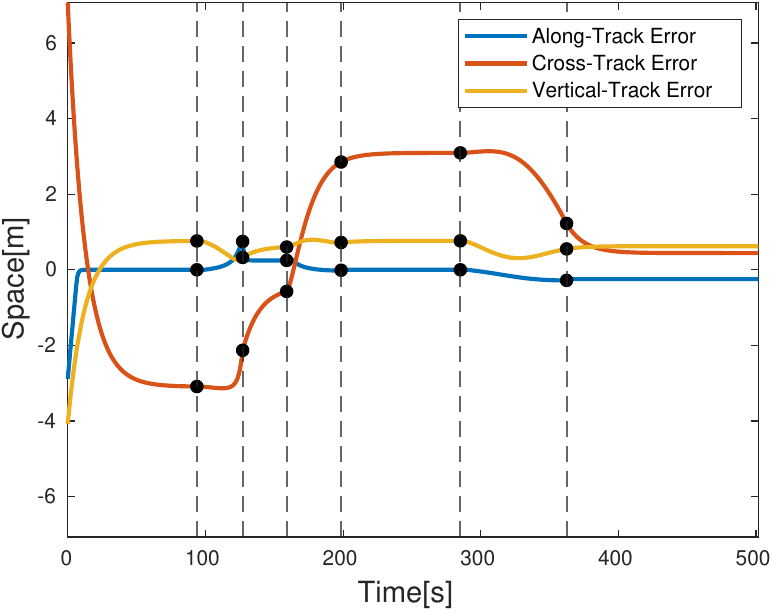}
\includegraphics[width=.315\linewidth]{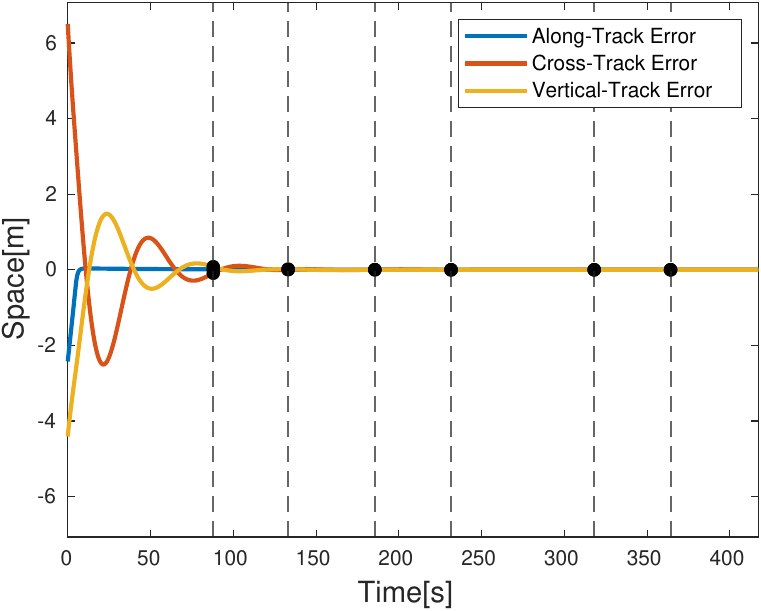}
%\centerline{\hspace*{.25cm} 
%(a) GL without current: track errors \hspace*{.7cm} 
%(b) GL with current: track errors \hspace*{.7cm} 
%(c) EGL with current: track errors}\bigskip\\
%\includegraphics[width=.315\linewidth]{img/base_segments_no_current/speed.pdf}
%%%\subfigure[GL with current: vehicle speed]{           
%\includegraphics[width=.315\linewidth]{img/base_segments_current/speed.pdf}
%\includegraphics[width=.315\linewidth]{img/extended_segment/speed.pdf}
%\smallskip\\
%\centerline{\hspace*{.2cm} 
%(d) GL without current: vehicle speed \hspace*{.1cm} 
%(e) GL with current: vehicle speed \hspace*{.1cm} 
%(f) EGL with current: vehicle speed}\bigskip%\\
\caption{Kinematic simulations: track errors obtained with the guidance laws of Section~\ref{sec:03} (GL) and \ref{sec:04} (EGL) for the first simulation scenario. The results are shown for a simulation without current (left) and with current when the basic guidance law (center) or its extension (right) are considered.}
    \label{fig:exm01-02}
\end{figure}

\begin{figure}[!t]
    \centering
%\subfigure[first simulation scenario]
{\includegraphics[width=.315\linewidth]{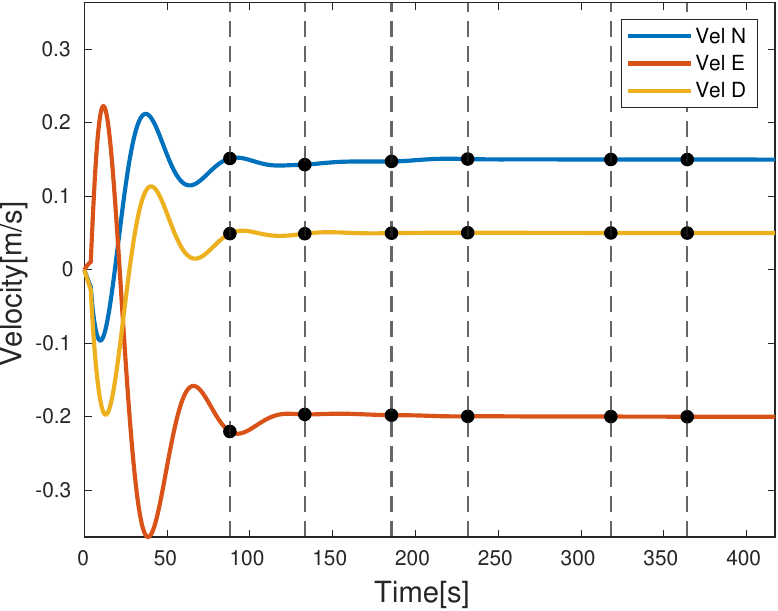}}
%\subfigure[second simulation scenario]
{\includegraphics[width=.315\linewidth]{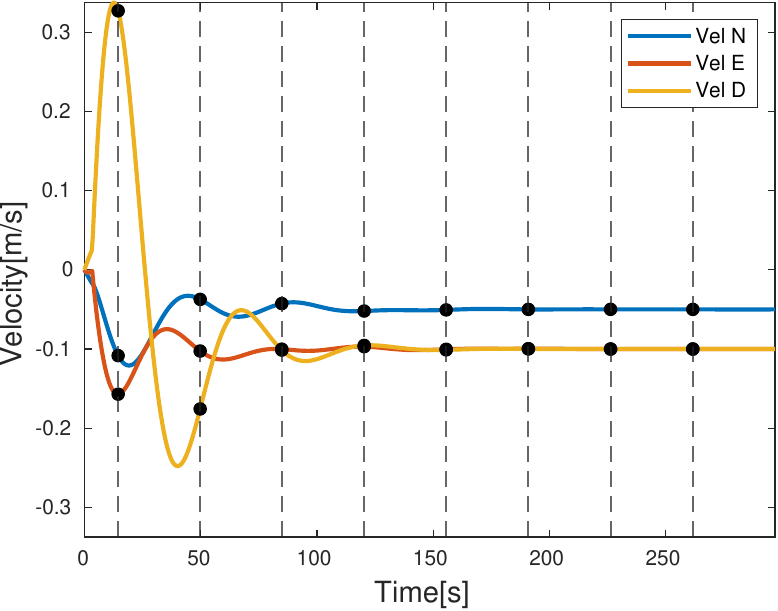}}
%\subfigure[third simulation scenario]
{\includegraphics[width=.315\linewidth]{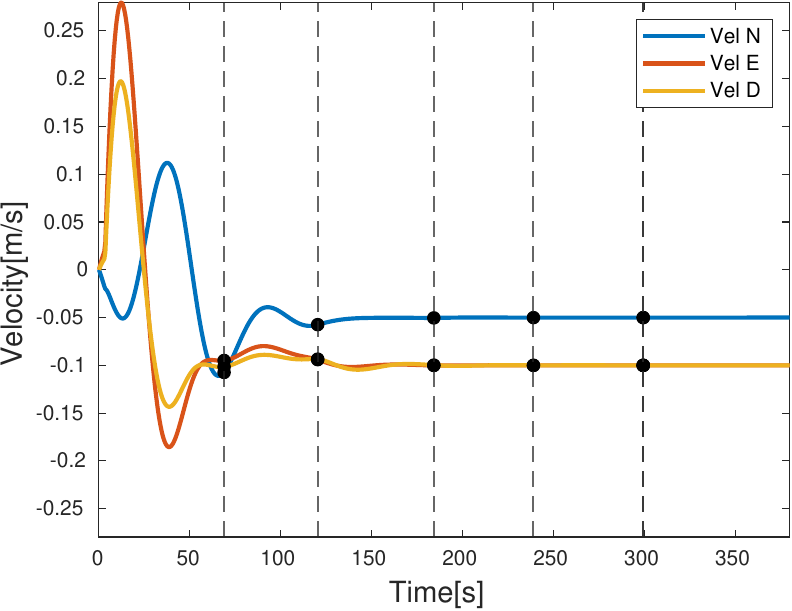}}
    \caption{Kinematic simulations: estimation of the current velocity in the first (left), second (center), and third (right) simulation scenario.
    %$\dot{\bm{\eta}}_c = (0.15 , -0.2 , 0.05)^\top \, \rm{m} / \rm{s}$. %The black dots and the dotted vertical lines indicate the interpolation points, in Nord-East-Down reference frame
    }
\label{fig:exm01-03}
\end{figure}

%\subsection{First simulation scenario}
%\label{section:segmentSettings}
In the first scenario, the eight points
$
\bm{p}_{0} = (0, 0, 0)^\top, \bm{p}_{1} = (20, 20, 20)^\top, \bm{p}_{2} = (35, 25, 25)^\top, \bm{p}_{3} = (50, 10, 25)^\top, 
\bm{p}_{4} = (45, -5, 20)^\top, \bm{p}_{5} = (25, -25, 0)^\top, \bm{p}_{6} = (10, -30, -5)^\top, \bm{p}_{7} = -(5,15,5)^\top
$
are interpolated for the construction of a  $C^1$ PH spline. The tangents in these points are choosen as $\bm{d}_j = \bm{d}_{j+1} =(\bm{p}_{j+1} - \bm{p}_{j})/\Vert\bm{p}_{j+1} - \bm{p}_{j}\Vert, j = 0,2,4,6$. The result is a $C^1$ PH spline curve composed by seven segments, with total length $L = 167.0983\, \rm{m}$. By considering as starting point of the vehicle $\bm{\eta}(0)=(-5, 5, -5)^\top \rm{m}$, we test the guidance law (GL) and its extended version (EGL). The drifting current is $\dot{\bm{\eta}}_c = (0.15 , -0.2 , 0.05)^\top \, \rm{m}/\rm{s}$. The results are shown in Fig.~\ref{fig:exm01-01}. The clear shift between the path and the trajectory of the vehicle (center plot in Fig.~\ref{fig:exm01-01}), exactly in the direction of the current vector, demonstrates that the basic GL is not able to handle current disturbances. This issue is resolved when the EGL is considered (right plot in Fig.~\ref{fig:exm01-01}). Fig.~\ref{fig:exm01-02} (center) confirms that the error does not converge when the GL is considered in presence of current, since its value depends on the relative direction of the path with respect to the current drift.
%. In the two segments where the tangent direction is nearly aligned to the current vector, the errors are smaller, while in the two segments where the tangent direction is perpendicular to the current vector, the errors grow. 
%Fig.~\ref{fig:exm01-02}(e) shows that also the velocity of the vehicle is highly affected by the presence of the current and depends on the direction of the path. 
%As expected, in the third and last spline segment, whose tangent vectors have opposite directions, the value of the speed is higher or lower if the vehicle moves in the direction of the current or not. 
%As a result, even if the two segments have the same length, the time required to cover them is totally different.  
%The starting point of the vehicle and the drifting current are the same as before. 
In Fig.~\ref{fig:exm01-02} (right) it is also possible to see that, when the EGL is considered, the error converges to zero with just some overshot effects, as normal consequences of the integral action in the current estimator update. Finally Fig.~\ref{fig:exm01-03} (left) shows that the current estimation converges to the value of current introduced in the kinematics.
\begin{figure}[t]
\includegraphics[width=.325\linewidth]{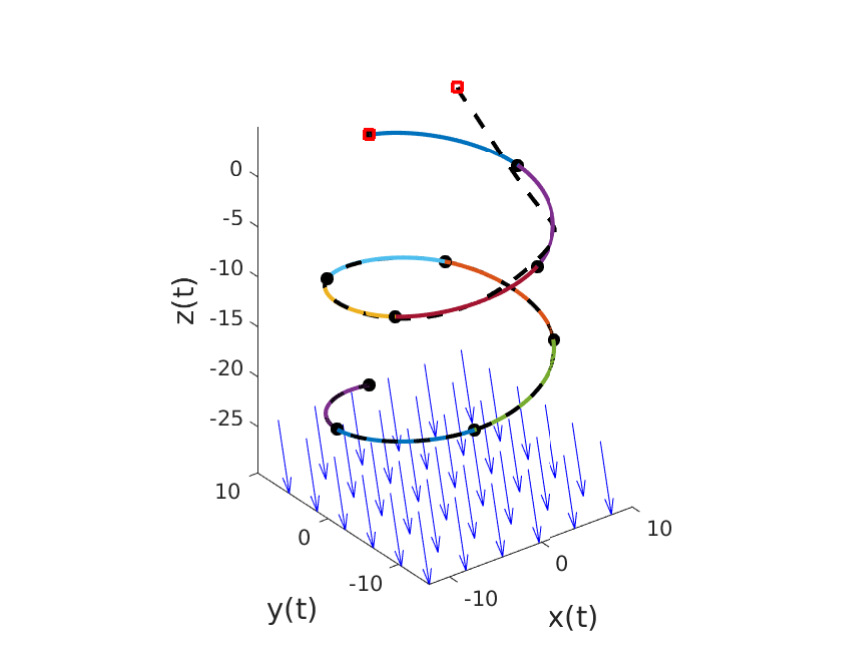}
{\includegraphics[width=.325\linewidth]{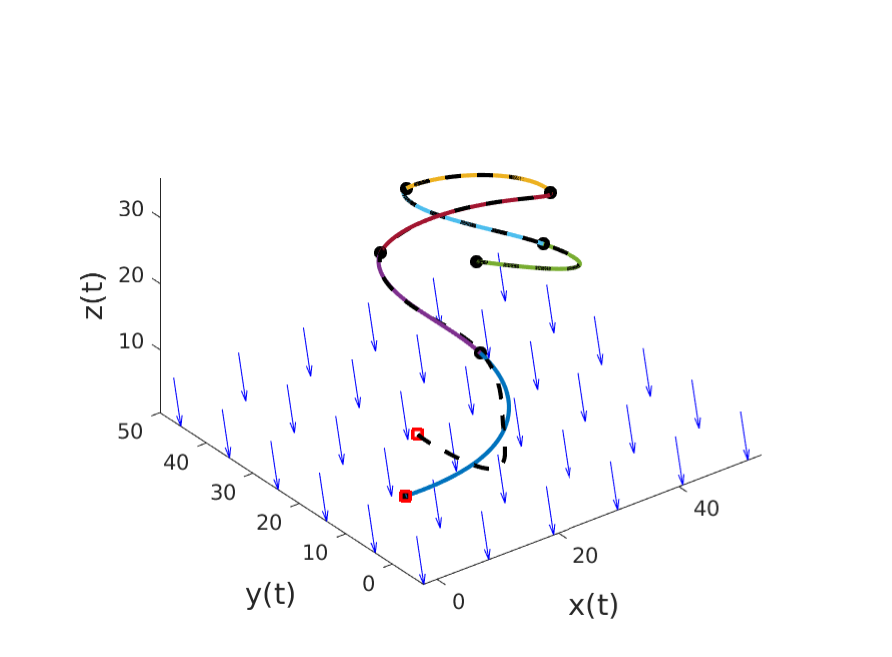}}
\includegraphics[width=.325\linewidth]{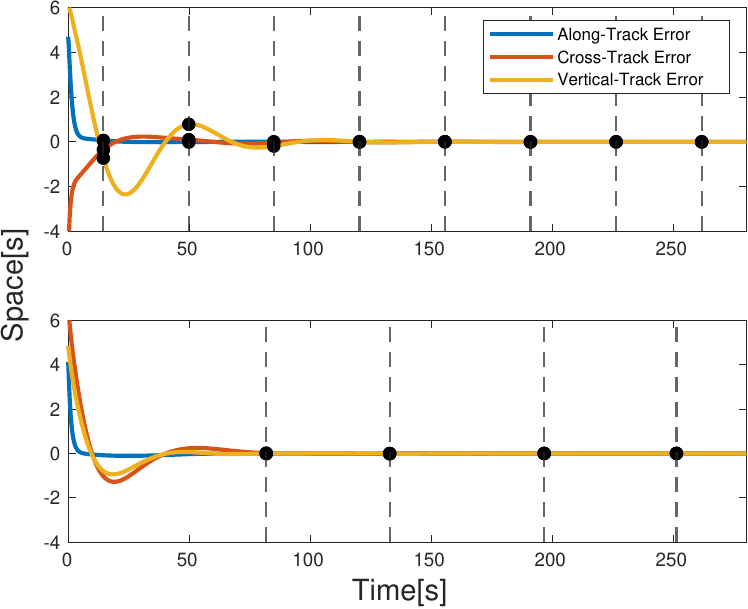}
\caption{Kinematic simulations: results obtained considering the extended guidance laws of Section~\ref{sec:04} for the last two scenarios. On the left (second scenario) and on the center (third scenario), the paths to be followed (solid line) are shown together with the paths of the vehicle (dashed line). The track errors obtained for the second (top) and third (bottom) simulation scenario are also shown (right).}
\label{fig:exm02and03}
\end{figure}

For the second scenario, we sample ten points and derivative vectors from a 
circular helix of the form 
$
\mathbf{x}(s) = 
\left(a\sin({s}/{s_0}),\,
a\cos({s}/{s_0}),\,
-\,b({s}/{s_0})
\right)^\top,
$
with 
$
s_0 = \sqrt{(a^2+b^2)},\, a=10, \, b=2,
$ 
at equidistant parameter values 
$
s \in [0 \,, \, n(2 \pi s_0)\,]\,,
$ with $n=2$ number of turns. The resulting $C^1$ PH spline curve with total length $L = 127.2619 \,\rm{m}$, computed with the CC selection criterion, is shown in Fig.~\ref{fig:exm02and03} (left) for a vehicle starting point equals to $\bm{\eta}(0) = (5,5,5)^\top \, \rm{m}$. The result obtained with the extended guidance law is presented in the same figure, considering a current velocity $\dot{\bm{\eta}}_c = (-0.05,-0.1,-0.1) \, \rm{m} / \rm{s}$. Again, the EGL is able to preserve the robustness of the scheme and allows to follow the desired path, as confirmed by the track errors in Fig.~\ref{fig:exm02and03} (right, top plot) and by the current estimation in Fig.~\ref{fig:exm01-03} (center).

In the third scenario, we construct the $C^1$ PH spline interpolating the following seven points: 
$
\bm{p}_{0} = (0, 0, 10)^\top,  
\bm{p}_{1} = (20, 10, 20)^\top, 
\bm{p}_{2} = (15, 25, 30)^\top,$
$
\bm{p}_{3} = (37, 17, 35)^\top, 
\bm{p}_{4} = (27, 35, 31)^\top, 
\bm{p}_{5} = (45, 29, 19)^\top, 
\bm{p}_{6} = (50, 50, 5)^\top.
$
The tangents in these points are chosen equal to the tangents  obtained with the $C^2$ cubic spline interpolant. The result is a $C^1$ PH spline curve composed by six segments, with total length $L = 153.6724 \, \rm{m}$. By considering $\bm{\eta}(0)=(5, 5, 15)^\top \, \rm{m}$ as starting point of the vehicle, we test the extended version of the guidance law with a drifting current $\dot{\bm{\eta}}_c = (-0.05,-0.1,-0.1) \, \rm{m}/\rm{s}$, see Fig.~\ref{fig:exm02and03} (center). Fig.~\ref{fig:exm02and03} (left, bottom plot) shows that the track-errors converge to zero while Fig.~\ref{fig:exm01-03} (right) confirms that the current estimation converges to the real current.
In table \ref{tab:estimatedTimeOfArrival} we compare the simulated time of arrival of the third simulation scenario with the estimated one. The latter is computed as $L_{r}(u)/U_0$, where $L_{r}(u)$ is the length of the portion of the trajectory from the point $\eta_p(u)$ to the end of the curve. It is clear that, as the vehicle approaches the trajectory, the estimate gets closer to the real data. Finally, we compare the computation time of the arc length using the exact formula available for PH splines with its computation time using the \texttt{integral} function of MATLAB. The result confirms that the first computation is 10-11 times faster than the second one. 

%\begin{figure}[t]
%\centering
%{\includegraphics[width=.325\linewidth]{img/third_scenario/path_following.pdf}}
%{\includegraphics[width=.325\linewidth]{img/third_scenario/lin_vel.pdf}}
%{\includegraphics[width=.325\linewidth]{img/third_scenario/ang_vel.pdf}}
%\caption{Third simulation scenario: results obtained with the guidance law of Section~\ref{sec:04} (EGL). 
%The path to be followed (solid line) is shown together with the path of the vehicle (dashed line) in a simulation (left), together with its relative velocity (center) and angular velocity (right).}
%\label{fig:exm03-01}
%\end{figure}

\begin{table}[t]
\centering
%\begin{tabular}{ p{2.5cm} | p{1.5cm}| p{1.5cm} | p{1.5cm} | p{1.5cm} | p{1.5cm} |p{1.5cm}}
%\resizebox{\textwidth}{!}
{
\begin{tabular}{lllllll}%{lrrrrrr}%{c @{\extracolsep{\fill}} cccccc}
 \hline
Parameter value &0.01 & 22.72 & 45.43 & 68.15 & 90.87 & 113.58 \\
 \hline
 Real ToA  & 379.78 s & 316.16 s & 255.46 s & 194.66 s & 136.63 s & 72.17 s \\
\hline 
 Estimated ToA & 384.11 s & 314.12 s & 255.66 s & 194.67 s & 136.64 s & 72.17 s \\
 \hline
\end{tabular}}
\medskip
\caption{Comparison between estimated time of arrival (Estimated ToA) and simulated time of arrival (Real ToA) in the third scenario for different parameter values.}
\label{tab:estimatedTimeOfArrival}
\end{table}

\subsection{Dynamic simulatons}

We finally validates the performance of the extended guidance law, presenting the results of numerical simulations based on the dynamic model of a real AUV called ``Zeno'' \cite{gmmprma2018}. The behaviour of a generic marine vehicle is described by the following vectorial equation of 6 components (corresponding to the 6 DoFs),
\[
\bm{M}\dot{\bm{\nu}}+\bm{C}(\bm{\nu})\bm{\nu}+\bm{D}(\bm{\nu})\bm{\nu}+\bm{g}(\mathcal{Q})=\bm{\tau},
\]
together with equation \eqref{fossen_kinematics_pos} and \eqref{fossen_kinematics_ori}. For a comprehensive description of all the involved quantities see \cite{fossen1994}. 
In the specific case of Zeno, the numerical values of all the dynamic model parameters have been identified in \cite{sacco2020}. 
%and here reported in Table ~\ref{tab:ZenoModel}. The measurement units have been omitted in order to lighten the notation.
%\begin{table}[t]
%\centering
%%\resizebox{\columnwidth}{!}
%\begin{tabular}{lclc}%\begin{normalsize}
% \hline
% $m$ & $40.25$
% &$\bm{I}_0$ & $\rm{diag}(0.94,3.28,3.79)$\\
% \hline 
% $\bm{A}_{11}$ & $\rm{diag}(87.75,72.62,85.07)$
% &$\bm{A}_{12}$ & $\bm{0}_{3 \times 3}$\\
% \hline
% $\bm{A}_{21}$ & $\bm{0}_{3 \times 3}$
% &$\bm{A}_{22}$ & $\rm{diag}(2.30,5.16,5.08)$\\
% \hline
%  $G$ & $394.69$
%  &$B$ & $396.69$\\
%  \hline
%  $\bm{r}_B$ & $(0 , 0 , -0.0015)^\top$\\
%  \hline  
%  $\bm{D}_L$ & $\rm{diag}(25.27,57.99,75.73,1.58,1.81,1.15)$ \\
%  \hline
%  $\bm{D}_Q(\bm{\nu})$ & $\rm{diag}(35.77,38.12,0,0.78,6.45,9.66)\rm{diag}(\vert \bm{\nu}\vert)$\\
% \hline
%%\end{normalsize}
%\end{tabular}
%\caption{Numerical values of the parameters identifying the  dynamic model of Zeno.}
%\label{tab:ZenoModel}
%\end{table}
Any dynamic simulation requires a suitable controller, whose purpose is to compute the vector of thrusts $\bm{\tau}$ that the vehicle has to set up in order to follow the kinematic references provided by the guidance law. We choose to design the following model-based controller:
\begin{equation*}
\bm{\tau} = \bm{M}\bm{K}_p \bm{e}_p + \bm{C}(\bm{\nu})\bm{\nu} + \bm{D}(\bm{\nu}_d)\bm{\nu}_d,
\end{equation*}  
where $\bm{\nu}_d =\left((\bm{v}_{rd}^b)^\top , (\bm{\omega}_d^b)^\top\right)^\top$ is the vector of references provided by the extended guidance law introduced in Section \ref{sec:04}, $\bm{e}_p = \bm{\nu}_d - \bm{\nu}$ is the velocity error, and $\bm{K}_p = \rm{diag}(20 , 0 , 0 , 0 , 10 , 10)$ is a $6\times6$ positive definite diagonal matrix containing the tunable parameter for the proportional terms of the controller. Furthermore, we consider to deal with a generic under-actuated vehicle, providing to the simulator a required effort $\bm{\tau} = (\tau_x , 0 , 0 , 0 , \tau_p , \tau_r)^\top$ which neglects the terms of sway, heave, and roll. Finally, a saturation on the thrusts has been introduced, based on the electromechanical limits of the real vehicle, In particular, we set a limit of $80 \rm{N}$ for surge thrust and a limit of $10 \rm{N \cdot m}$ for pitch and yaw torque.  
 
\begin{figure}[t]
\centering
\includegraphics[width=.325\linewidth]{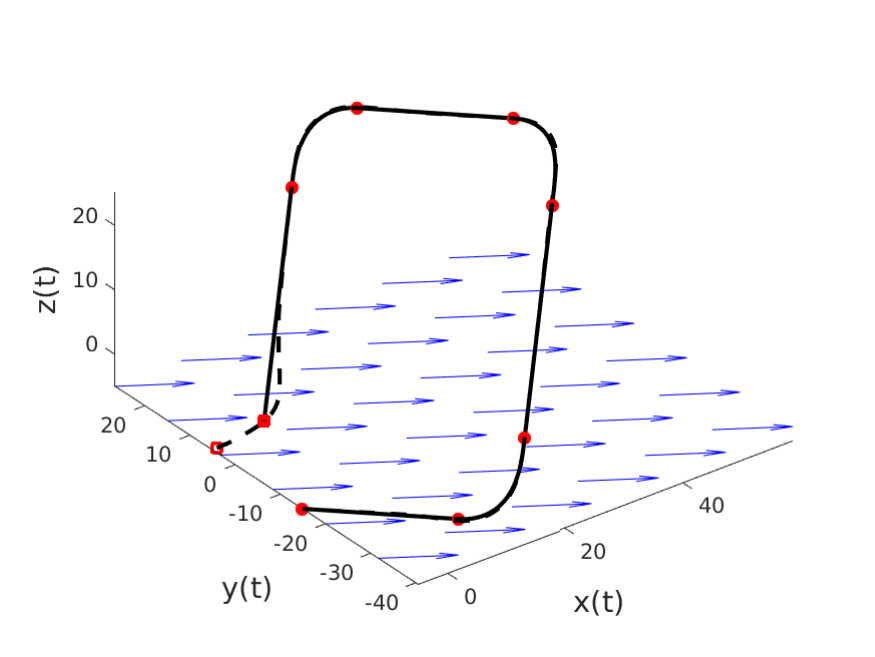}
\includegraphics[width=.325\linewidth]{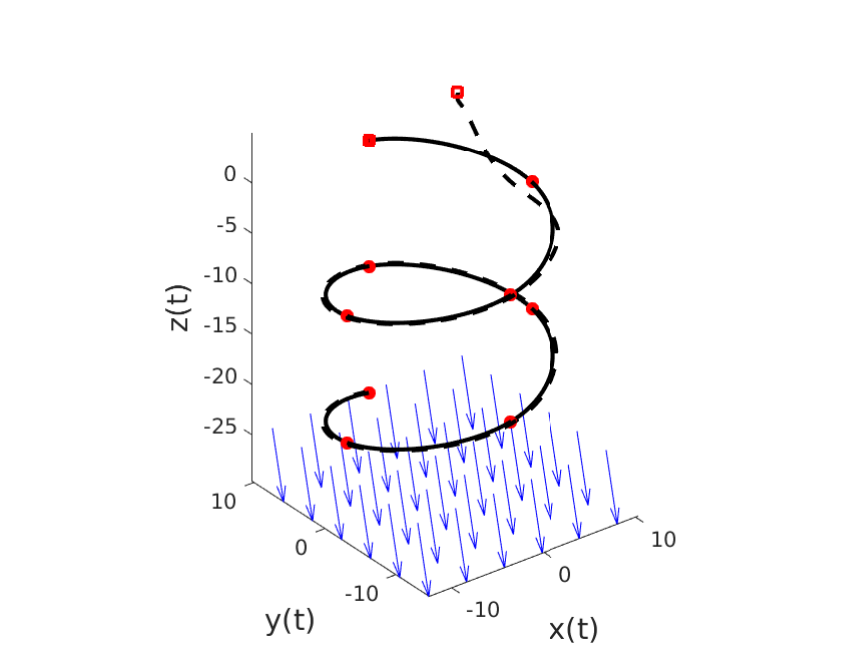}
\includegraphics[width=.325\linewidth]{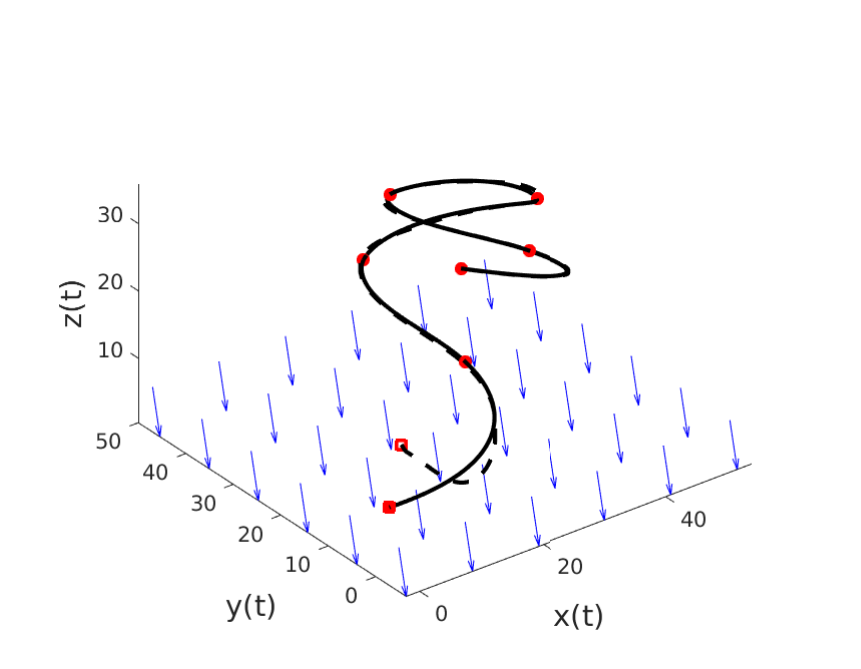}

\caption{Dynamic simulations: the path to be followed (solid line) is shown together with the path of the vehicle (dashed line) for the first (left), the second (center), and the third (right) simulation scenario.}
    \label{fig:dyn01}
\end{figure}

\begin{figure}[t]
\centering
\includegraphics[width=.325\linewidth]{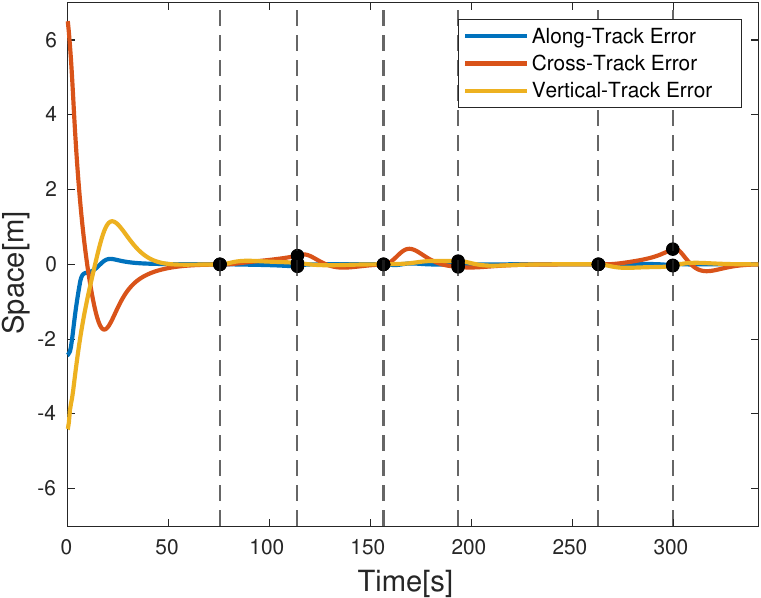}
\includegraphics[width=.325\linewidth]{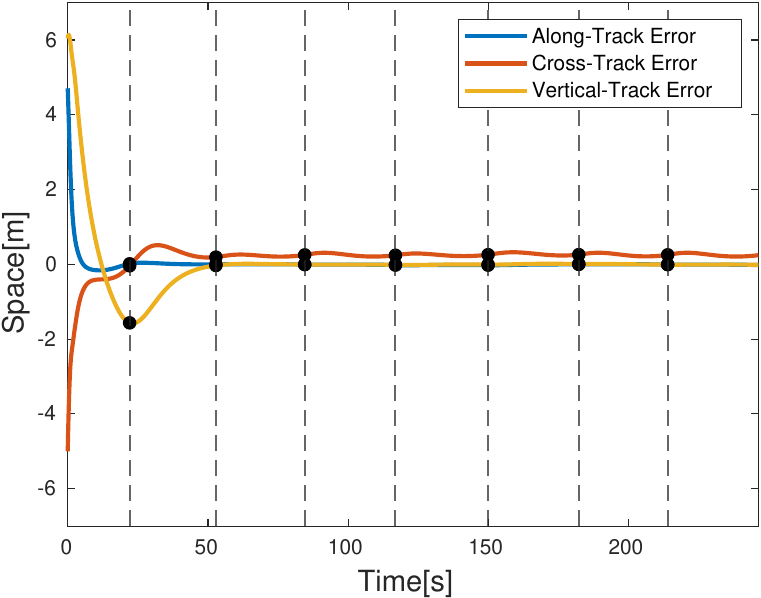}
\includegraphics[width=.325\linewidth]{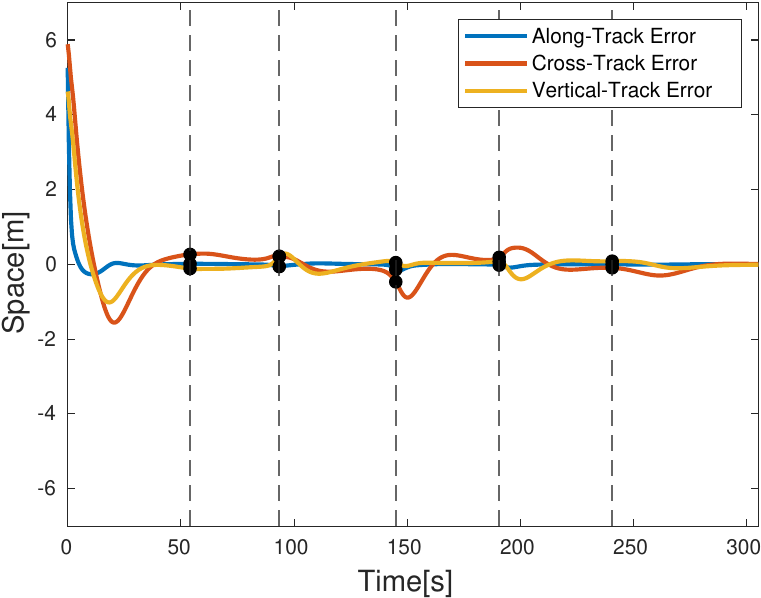}

\caption{Dynamic simulations: track errors for the first (left), the second (center), and the third (right) simulation scenario.}
    \label{fig:dyn02}
\end{figure}

We simulate the behavior of Zeno by considering all the three simulation scenarios. 
%The guidance parameter introduced in Section~\ref{sec:03} and~\ref{sec:04} are set as follows: $\gamma = 1 \,\rm{s}^{-1}$, $U_0 = 0.5 \, \rm{m}/\rm{s}$, $\Delta_e = 2.5 \, \rm{m}$ and $\mu = 1$, $\bm{O}^b_h = (0.8,0,0)^\top \, \rm{m}$ and $\kappa_c=0.015 \,\rm{s}^{-2}$. 
In Fig.~\ref{fig:dyn01} it is possible to verify the convergence of the vehicle on the path. The corresponding track errors are reported in Fig.~\ref{fig:dyn02}. It is remarkable to note that in some cases the cross track error grows again. This effect is connected to the electromechanical limits of the vehicle, which is not able to follow high curvature profiles. However, the error decreases as soon as the curvature gets lower.

\section{Closure}
\label{sec:07}
We presented a robust path following scheme for autonomous vehicles, which allows to handle disturbance effects caused by currents drifting, also suited for under-actuated vehicles. The theoretical validity of the scheme was proved through geometrical and error estimation convergence. Finally, we confirmed the results with simulations based on tangent continuous curvilinear paths obtained by Hermite interpolation with spatial $C^1$ PH quintic splines, which enable accurate and efficient arrival time estimations. 
We also tested the guidance law on a classic scenario, considering a 6 DoFs underwater vehicle subject to thrust limitation and decoupled dynamical speed controller. The Zeno AUV model, identified with experimental tests, was used to evaluate the combination of a standard existing controller with the new guidance law under sea current operating condition. The proposed guidance law eases the controller from any integral action and improves convergence by estimating the contribution of the sea current directly into the navigation frame.
\bibliographystyle{plain}
\bibliography{biblio}
\end{document}